\documentclass[12pt]{article}
\usepackage{geometry}                
\usepackage{graphicx}
\usepackage{amssymb}					
\usepackage{amsmath}					
\usepackage{mathrsfs}
\usepackage{tikz}
\usepackage{epstopdf}
\usepackage{amsthm}
\newtheorem*{theorem*}{Theorem}
\newtheorem*{lemma*}{Lemma}
\newtheorem*{example*}{Example}
\newtheorem{theorem}{Theorem}
\newtheorem{lemma}{Lemma}

\newtheorem{proposition}{Proposition}

\def\b0{{\bf 0}}
\def\b1{{\bf 1}}

\def\cC{{\cal C}}

\def\cP{{\cal P}}

\def\cE{{\cal E}}

\def\cC{{\cal C}}

\def\cG{{\cal G}}

\def\n{\noindent}

\begin{document}
\title{
Construction numbers: How to build a graph \thanks{to appear in {\it Missouri J. of Mathematical Sciences}} }
\author{Paul C. Kainen,  \tt kainen@georgetown.edu
}
\date{}                                           

\newcommand{\Addresses}{{
  \bigskip
  \footnotesize

\vspace{-0.27cm}

\n

\par\nopagebreak
}}

\maketitle
\begin{abstract}
\n A construction sequence for a graph is a listing of the elements of the graph (the set of vertices and edges) such that each edge follows both its endpoints.  The construction number of the graph is the number of such sequences.  We determine this number for various graph families.
\end{abstract}

\section{Introduction}

How many ways are there to enumerate the elements (vertices and edges) of a graph so that each edge is listed after its endpoints?  This {\it construction number} (c-number for short) turns out to be interesting but hard to compute.


In fact, c-number is the number of ways to extend the natural poset of a graph, given by its elements under inclusion, to a linear order.  
The set of linearizations of a poset were studied by Stanley
\cite{stanley}, but he 
only much later \cite{two} enumerated c-numbers of paths.  Here, we determine these numbers also for several other basic graph families: cycles, stars, double-stars, complete graphs and disjoint unions of edges.
We obtain instances of well-known integer sequences.


Construction number constitutes a measure of graph complexity and grows quite rapidly.  For example, for the star graph with $n$ peripheral vertices, having $2n+1$ elements, the construction number is $c(K_{1,n}) = 2^n (n!)^2$, so $c(K_{1,10}) >  1.348 \times  10^{16}$.

{\it How does the topology of a graph relate to its construction number}? 
For trees with a fixed number of elements, paths have the lowest and stars have the highest construction number \cite{random-constr}.
The corresponding problems are open for maximal planar and maximal outerplanar graphs.

Section 2 has definitions, an example, and basic lemmas, while Section 3 and the Appendix give the c-number calculations.  There is a discussion in Section 4. 



\section{Basic definitions and lemmas}

For any set $S$, $|S|$ denotes cardinality, always finite in this paper.
For $G=(V,E)$ a graph, write $p := p(G) := |V|$ for {\bf order} of $G$, $q := q(G) := |E|$ for {\bf size} of $G$. 
Maximum and minimum vertex-degree are denoted $\Delta$ and $\delta$.  Let $[k] := \{1,\ldots, k\}$, $k \geq 1$, and let $\mathbb{N}_{>0}$ be the set of positive integers. Write ``$*$'' for concatenation.

Let $P_p$ be the path with $p$ vertices and $K_{1,n}$ the star with a hub of degree $n$ each of whose neighbors is an endpoint so $p(K_{1,n}) = n{+}1$; $C_p$ and $K_p$ are the cycle graph and complete graph with $p$ vertices.  Graph terminology is as in \cite{harary}. Let $\ell = p+q$.

The {\bf elements} of a graph $G = (V,E)$ are the set $V  \cup E$ of vertices and edges.  If $X$ is any set of element of $G$, the {\bf subgraph $G(X)$ induced by $G$} is the graph whose edges are the edges in
X and whose vertices are the vertices in X union the endpoints of the edges in X.

A {\bf construction
sequence} (or {\bf c-sequence}) for $G$ is a listing of all the elements without repetitions such that {\it each edge appears after both of its endpoints}.   For instance, the path $P_3 := (\{1,2,3\},\{12, 23\})$ has  construction sequences $(1,2,3,23,12)$ and $(1,2,12,3,23)$, while $(1,3,12,2,23)$ is {\it not} a c-sequence. The full set of 16 c-sequences of $P_3$ are shown later in this section.

Let $\cC(G)$ be the set of all c-sequences for $G=(V,E)$.
The {\bf construction number} 
$$c(G) :=|\cC(G)|$$ is  the number of distinct construction sequences.  So $c(P_3)=16$.
If $x \in \cC(G)$ and
$1 \leq i \leq \ell$, we let $x_{(i)} := (x_1, \ldots, x_i)$ and $G_i := G(\{x_1, \ldots, x_i\})$.  

A useful subset $\cE(G)$ of $\cC(G)$ is formed by the {\bf easy} c-sequences for which all vertices are listed first, then all the edges.  Other construction sequences are called {\bf non-easy}.  A superset $\cP(G)$ of $\cC(G)$ is formed by the arbitrary lists of elements without repetitions but with no constraint on when edges may appear.
Hence, 
\begin{equation}
\cE(G) \subseteq \cC(G) \subseteq \cP(G).
\label{eq:sets-veb}
\end{equation}
and as an immediate consequence, for each graph $G$, 
\begin{equation}
p!q! \leq c(G) \leq (p+q)!
\label{eq:vertex-edge-bounds}
\end{equation}
We call these the {\bf vertex-edge bounds} on the construction number.

\begin{lemma} 
If $G$ is connected, then the last $\delta(G)$ entries in any construction sequence $x \in \cC(G)$ are edges.  Moreover, if $x_i \in V(G)$, then as each vertex must be followed by the edges with which it is incident,
$\deg(v_i) \leq \ell -i$.
\label{lm:delta}
\end{lemma}


If $G' \subseteq G$ and $x \in \cP(G)$, then we define $x\,|_{G'}$, the {\bf restriction} of $x$ to $G'$, to be the element in $\cP(G')$ consisting of the entries in $x$ which are elements of $G'$.  When $x \in \cC(G)$, then $x\,|_{G'}$ is also a c-sequence.  What about the converse?

\begin{lemma}
Let $G' \subseteq G$. 
If $y \in \cP(G)$ with $y|_{G'} \in \cC(G')$, then $y \in \cC(G)$ if and only if for each edge $e=uw$ of $G$ that is not in $G'$, $e$ follows $u$ and $w$ in the listing $y$. 
\end{lemma}

The number of ways to extend a c-sequence for $G$ - $v$ to a c-sequence for $G$  depends on the particular construction sequence.  
For example, take $P_2 = (\{1,2\},\{a\})$, where $a$ is the edge $12$. Then $\cC(P_2)=\{x',y'\}$, where $x':=(1,2,a) \equiv 12a$ (in {\it short form}) and $y' \equiv 21a$.  Consider $P_3 = (\{1,2,3\}, \{a,b\})$, where $b=23$. As $P_2 \subset P_3$, each c-sequence for $P_3$ extends a c-sequence of $P_2$.  One finds that
$x'$ has exactly 7 extensions to an element in $\cC(P_3)$ (only one of which is non-easy)
\[ 312ab, 312ba,132ab,132ba,123ab,123ba,12a3b,\]
while $y'$ has exactly 9 extensions (two of which are non-easy)
\[321ab,321ba,32b1a,231ab,231ba,23b1a,213ab,213ba,21a3b.   \]
These are the 16 elements of $\cC(P_3)$.

Given two elementwise disjoint finite sequences $s_1$ and $s_2$ of lengths $n$ and $m$, we define a {\bf shuffle} of the two sequences to be a sequence of length $n+m$ which contains both $s_1$ and $s_2$ as subsequences. The number of shuffle sequences of $s_1$ and $s_2$ is ${{n+m}\choose{n}}$, giving the construction number of a disjoint union in terms of its parts.
\begin{lemma}
If $x_1$ and $x_2$ are c-sequences for disjoint graphs $G_1$ and $G_2$, resp., then each shuffle of $x_1$ and $x_2$ is a c-sequence for $G_1 \cup G_2$, and conversely; hence,
\begin{equation}
c(G_1 \cup G_2) = c(G_1) c(G_2) {{\ell_1+\ell_2}\choose{\ell_1}},
\label{eq:disj-form}
\end{equation}
where $\ell_1$ and $\ell_2$ are the lengths of the sequences $x_1$ and $x_2$, resp.
\label{lm:union}
\end{lemma}

The previous lemma extends to any finite pairwise-disjoint union of graphs. If $G_1, \ldots, G_n$ have $\ell_1, \ldots, \ell_n$ elements, then for $\ell := \ell_1 + \cdots + \ell_n$, replacing binomial by multinomial coeffient, the pairwise-disjoint union ({\it coproduct}) $G := \coprod_{i=1}^n G_i$ satisfies
\begin{equation}
c(G) = \prod_{i=1}^n c(G_i) {{\ell}\choose{\ell_1, \cdots, \ell_n}}.
\label{eq:disj-union}
\end{equation}

Construction numbers are equal for isomorphic graphs.  More exactly, we have
\begin{lemma}
Let $\phi: G \to H$ be an isomorphism of graphs. Then there is an induced bijection $\hat{\phi}: \cC(G) \to \cC(H)$ given by $\hat{\phi}(x) = \hat{x}$, where $\hat{x}_j = \phi(x_j)$ for $1 \leq j \leq \ell(G)$.
\label{lm:bij}
\end{lemma}


%
 
Let 
 $$\cC(v,G) := \{x \in \cC(G): x_1 = v\} \;\;\mbox{and}\;\; c(v,G) := |\cC(G,v)|.$$
%
%
%
%
%
%

A graph $G$ is {\bf vertex-transitive} (or {\bf edge-transitive})
if for every pair $(v,w)$ of vertices (or $(e,f)$ of edges), there is an automorphism carrying $v$ to $w$ (or $e$ to $f$).
\begin{proposition}
If $G = (V,E)$ is vertex transitive, then $c(G) = |V| \cdot c(v,G)$.
\end{proposition}

%
%
For $e \in E$ 
write
$\cC(G,e) := \{x \in \cC(G): x_\ell = e\}$; put $c(G,e) := |\cC(G,e)|$.
Every construction sequence for a nontrivial connected graph ends in an edge. So we have
\begin{lemma}
Let $G=(V,E)$ be nontrivial connected. Then 
\begin{equation}
c(G) = \sum_{e \in E} c(G - e).
\label{eq:edge-recur}
\end{equation}
\end{lemma}
\begin{proof}
Each $y' \in \cC(G,e)$ has the form $y' = y * e$ for $y \in \cC(G - e)$.
\end{proof}
%
\begin{proposition}
Let $G=(V,E)$ be nontrivial connected and edge transitive. Then $c(G) = |E| \cdot c(G,e)$.
\label{pr:e-trans}
\end{proposition}

In fact, for a connected graph, the number of edges concluding any construction sequence is at least the minimum degree of the graph.

\section{Construction numbers for some graph families}

Construction numbers are determined for stars and double-stars, for paths and cycles, and for disjoint unions of edges. For complete graphs, see the Appendix. 
\begin{theorem}[c-number of a star]
For $n \geq 0$, $c(K_{1,n}) = 2^n(n!)^2$.
\end{theorem}
\begin{proof}
For $n=0,1$, the result holds. Suppose $n \geq 2$ and let $x = (x_1, \ldots, x_{2n+1})$ be a construction sequence for $K_{1,n}$. There are $n$ edges $e_i = v_0 v_i$, where $v_0$ is the central node, and one of the edges, say, $e_i$, must be the last term in $x$.  This leaves $2n$ coordinates in $x' := (x_1, \ldots, x_{2n})$ and one of them is $v_i$.  The remaining $(2n-1)$ coordinates are a construction sequence for the $(n-1)$-star $K_{1,n-1}$.  Hence, $c(K_{1,n})  = n (2n) c(K_{1,n} - v_i)= 2n^2 2^{n-1}(n-1)!^2 = 2^n (n!)^2$ by induction.
\end{proof}
The numbers  2, 16, 288, 9216, 460800 generated by the above formula count the number of c-sequences for $K_{1,n}$ for $n \in \{1,2,3,4,5\}$.
These numbers are the absolute value of the sequence \cite[A055546]{oeis}; they describe the number of ways to seat $n$ cats and $n$ dogs in a roller coaster with $n$ rows, where each row has two seats which must be occupied by a cat and a dog.

Interestingly, the value of $c(K_{1,n})$ 
asymptotically agrees with
the geometric mean of the vertex-edge bounds (\ref{eq:vertex-edge-bounds}) up to a constant multiple of $n^{3/2}$.  For we have 
\begin{equation}
 \frac{c(K_{1,n})}{(n+1)! n!} =\frac{2^n (n!)^2}{(n+1)!n!} = \frac{2^n}{n}
\end{equation}
while
\begin{equation}
 \frac{(2n+1)!}{c(K_{1,n})} = \frac{(2n+1)!}{2^n(n!)^2} \sim \frac{2^{1+n} \sqrt{n}}{\sqrt{\pi}} = c \;2^n \sqrt{n}, \;\mbox{for}\;c= \frac{2}{\sqrt{\pi}} \approx 1.128....
\end{equation}

A {\bf double star} is a tree of diameter at most 3. For $a,b \geq 0$, let $D_{a,b}$ denote the double star formed by the union of $K_{1,a+1}$ and $K_{1,b+1}$ with one common edge.  So $D_{a,b}$ consists of two adjacent vertices, say $v$ and $w$, where $v$ is adjacent to $a$  and $w$ to $b$ additional vertices, and all vertices other than $v$ and $w$ have degree 1.
One gets a recursion for the construction number of a double star.

\begin{theorem}
Let $a,b \geq 0$ be integers. Then 
$$c(D_{a,b}) = f(a,b) + \beta(a,b) c(D_{a-1,b}) + \beta(b,a) c(D_{a,b-1}),$$
where $f(a,b) = c(K_{1,a} \coprod K_{1,b}) ={{2a + 2b + 2}\choose{2 a + 1}} 2^a (a!)^2 2^b (b !)^2$, and $\beta(a,b): = a(2a+2b+2)$.
\end{theorem}
We don't yet have a closed form but can calculate the values. For example, $D_{2,2}$, the $6$ vertex tree with two adjacent cubic nodes, has 402432 c-sequences. But $c(K_{1,5}) = 460800$.  

Stars maximize the construction number over all trees of a fixed order.  In \cite{random-constr}, it is shown further that
{\it the extremal trees have extremal construction numbers}.
\begin{lemma}
If $T$ is any $p$-vertex tree, then $c(P_p) \leq c(T) \leq c(K_{1,p-1})$.  
\end{lemma} 

For the path, we show below that $c(P_p) = T_p$, where $T_p$ is the $p$-th {\bf Tangent number} \cite[A000182]{oeis}.  The following formula 
\cite{mathw}, \cite[24.15.4]{dlmf} expresses Tangent number in terms of Bernoulli number
\begin{equation}
T_p = (1/p) {{2^{2p}}\choose{2}} |B_{2p}|,
\label{eq:path}
\end{equation}
where $B_{2p}$ is the $2p$-th Bernoulli number \cite[(A027641$/$A027642)]{oeis},

An asymptotic analysis shows  $c(P_{p})$ is exponentially small compared to $c(K_{1,p-1})$,
\begin{equation}
\frac{c(P_{p})}{c(K_{1,p-1})} \sim \Big(4e^{-2} \ \pi^{-1/2} \Big)  p^{1/2} \Big(\frac{8}{\pi^2}\Big)^p.
\label{eq:path-to-star}
\end{equation}

\smallskip

By Proposition \ref{pr:e-trans}, counting c-sequences for $C_p$ and $P_p$ are equivalent problems.
\begin{lemma}
If $p \geq 3$, then $c(C_p) = p \cdot c(P_p)$.
\label{lm:cycle-path}
\end{lemma}

Using Lemma \ref{lm:cycle-path} and equation (\ref{eq:path}), we have for $p \geq 3$,
\begin{equation}
c(C_p) = {{2^{2p}}\choose{2}} |B_{2p}|.
\label{eq:cycle}
\end{equation}
In fact, the formula holds for $p \geq 1$, where cycles are CW-complexes but not graphs: $C_1$ has one vertex and a  loop, and $C_2$ has two vertices
and two parallel edges.

Before determining $c(P_p)$, we give a Catalan-like recursion for these numbers.
\begin{lemma}
For $n \geq 2$, $c(P_n) = \sum_{k=1}^{n-1} c(P_k) \,c(P_{n-k}) \,{{2n-2}\choose{2k-1}}.$
\label{lm:cat}
\end{lemma}
\begin{proof}
Any c-sequence $x$ for $P_n$ has last entry an edge $e$, whose removal creates disjoint subpaths $P, P'$ with $k$ and $n-k$ vertices, resp., for some $k$, $1 \leq k \leq n-1$, and $x_{(2n-2)}$ is a c-sequence for $P \cup P'$. Now use Lemma \ref{lm:union}.
\end{proof}
Trivially, $c(P_1) = 1$.  By the recursion, one gets
$ 1, 2, 16, 272, 7936, 353792$ for $c(P_n)$, $n=1,\ldots,6$, empirically \cite[A000182]{oeis} the sequence of tangent numbers, $T_n$.  The following theorem was first shown by Stanley \cite{two}; we give two different proofs.


\begin{theorem}
For all $n \geq 1$, $c(P_n) = T_n$.
\end{theorem}
\begin{proof}
Let $U(n)$ denote the set of {\bf up-down} permutations
of $[n]$, 
where successive differences switch sign, and the first is positive.  According to \cite{mathw2}, $|U(2n-1)| = T_n$.

Dan  Ullman observed that there is a bijection from $\cC(P_n)$ to $U(2n-1)$.  Indeed, a permutation $\pi$ of the consecutively labeled elements of a path  is a construction sequence if and only if $\pi^{-1}$ is an up-down sequence.  Indeed, for $j = 1, \ldots, n-1$, $\pi^{-1}(2j)$ is the position in $\pi$ occupied by the $j$-th edge, while $\pi^{-1}(2j-1),\pi^{-1}(2j+1)$ correspond to the positions of the two vertices flanking the $2j$-th edge and so are smaller iff $\pi$ is a construction sequence. See also \cite[pp. 157 \& 180]{enum-comb1} \end{proof}
For example, 
labeling the elements of $P_5$ from left to right $(1,2,3,4,5,6,7,8,9)$, odd-numbers correspond to vertices and even numbers to edges, 
and $P_5$ has the c-sequence $\pi = (1,3,5,7,9,2,4,6,8)$, while
 $\pi^{-1} = (1,6,2,7,3,8,4,9,5)$
is up-down.)  

\smallskip

Here is a second proof that $c(P_n) = T_n$.
Let $J_r$ be the number of permutations of $\{0,1, \ldots, r+1\}$ which begin with `$1$', end with `$0$', and have consecutive differences which alternate in sign. Observe that $J_{2k}=0$ for $k \geq 1$ as the sequences counted by $J$ must begin with an {\it up} and end with a {\it down} and hence have an odd number of terms. These {\bf tremolo} sequences were shown by Street \cite{rs} to be in one-to-one correspondence with ``Joyce trees'' 
and to satisfy the following recursion.

\begin{proposition} [R. Street]
For $r \geq 3$,
$J_r = \sum_{m=0}^{r-1} {{r-1}\choose{m}} J_m J_{r-1-m}$.
\label{pr:street}
\end{proposition}
Street then proved that $J_{2n-1} = T_n$.
Hence, it suffices to show that $c(P_n) = J_{2n-1}$. Indeed,
$J_1 = c(P_1)$ and $J_3 = c(P_2)$. If you replace $J_{2r-1}$ by $c(P_r)$ and 
$J_{2r}$ by zero and re-index, then
the recursion in Proposition 5 becomes that of Lemma \ref{lm:cat}, so $c(P_n)$ and $J_{2n-1}$ both satisfy the same recursion and initial conditions and so are equal.\\

If $v$ is one of the endpoints of $P_n$, we can calculate $c(v, P_n)$ for the first few values, getting (with some care for the last term) $1,1,5,61$ for $n=1,2,3,4$.  In fact,
\begin{equation}
c(v, P_n) = S_n,
\label{eq:star-base}
\end{equation}
where $S_n$ is the $n$-th Secant number (\cite[A000364]{oeis}), counting the ``zig'' permutations.

The reader may try the following problem \cite{monthly-prob}: Find $c(K_n)$, where $K_n$ is the complete graph on $n$ vertices, before looking at its solution in the Appendix. 

Let $n K_2$ be the disjoint union of $n$ edges. By equation (\ref{eq:disj-union}), we have the following.
\begin{theorem} For $n \geq 1$,
$$c(n K_2) = 2^n {{3n}\choose{3,3,\dots,3}}.$$
\end{theorem}
The corresponding sequence is A210277 in \cite{oeis}.

\section{Discussion}

Stanley studied the number of linear extensions of a partial order \cite[p 8]{stanley}, using them to define a polynomial \cite[p 130]{enum-comb1}. In \cite[p 7]{two}, \cite[pp 157 \& 180]{enum-comb1} he showed the number of linear extensions of the partial order determined by a path is an Euler number, implying our results (\ref{eq:path}), (\ref{eq:cycle}), and (\ref{eq:star-base}) above. 

Brightwell \& Winkler \cite{bw91} showed that counting the linear extensions of a poset is $\#P$-complete and contrast this with randomized polynomial-time algorithms which estimate this number.  Their conjecture that $\#P$-completeness holds even for height-2 posets was proved by Dittmer \& Pak \cite{dp2020}.
Applications of linear extensions of posets to  equidistributed classes of permutations were given by Bj\"orner \& Wachs \cite{bw91}, and Burrow \cite{burrow} has studied using traversals of posets representing taxonomies and concept lattices to construct algorithms for information databases.

The set $[n]$ of positive integers under divisibility form a poset and the number $a(n)$ of linear extensions is enumerated in A016021 of \cite{oeis}. We note that understanding this integer sequence would be quite interesting as
\begin{equation}
a(n+1) = n a(n) \iff n+1 \;{\mbox{\it is prime}}.
\end{equation} 

A different approach to graph construction was motivated by the goal to describe the self-assembly of macromolecules performed by virus capsids in the host cell.  An {\it assembly tree} (Vince and Bon\'a \cite{vb2012}) is a gadget that builds up a graph from subgraphs induced by various subsets of the vertices. The number of assembly trees of a graph is the {\it assembly number}. See \cite{mw-an, mw-cn}.

An ``opportunity cost" can be defined by charging each edge for the number of steps in the c-sequence by which it follows the placements of both its endpoints.  One can consider min- and max-cost c-sequences for a given graph and algorithms that allows the graph to be built according to various rubrics. See \cite{gk}, with J. Gao, where we showed that max cost depends only on the degree sequence and that max and min cost have a ratio of three for the complete graph.

Graph construction constrained by minimizing the number of connected components and simplifying the graph information is studied in \cite{UK} with P. Ulrickson.  Further, we have used analytic methods in \cite{random-constr} with R. Stong.




%
It would be interesting to have actual construction sequences for biochemical and  neuronal networks, or the order in which atoms and bonds appear when  molecules are built. 


\subsection*{Acknowledgement}
The author thanks Stan Wagon for helpful suggestions and for noting reference \cite{vb2012}. He also thanks Richard Hammack for pointing out equation (\ref{eq:edge-recur}) and he is grateful to the referee for substantially improving the readability.

\subsection*{Appendix: Solution of the Monthly problem \cite{monthly-prob}}

\begin{theorem}
For $n \geq 1$ and letting $(2n-1)!! := (2n-1)(2n-3) \cdots  3 \cdot 1$, we have\\
\begin{equation}
c(K_n) = {{n+1}\choose{2}}!\Big/(2n-1)!!\;.
\end{equation}
\label{th:Kn}
\end{theorem}
Thus, $(c(K_n))_{n=1}^\infty$ is rapidly growing:
$1, 2, 48, 34560, 1383782400, ...\;$.
The first four terms agree with a sequence in the OEIS, but the fifth and subsequent terms do not.   As $\ell(K_n) = {{n+1}\choose{2}}$ 
the probability that a uniformly randomly selected sequence of vertices and edges is a c-sequence for $K_n$ is $1/(2n-1)!!$.  
The following proof of Theorem \ref{th:Kn} is due to Richard Stong. 
A different argument is given by Tauraso \cite{tauraso}.
\begin{proof}
Index the vertices by $1, \dots, n$.
Decompose the elements of $K_n$ into $n$ {\bf teams} $T_1, \ldots, T_n$, where the team $T_k := \{k\} \,\cup\,
\{ik \,:\,1 \leq i \leq k\}$
consists of vertex $k$ and the edges to it from vertices of lower index. 
So $T_1 = \{1\}$, $T_2 = \{2, 12\}$, etc.  For $1 \leq k \leq n$, team $T_k$ has $k$ members.  
The $n$ teams are pairwise-disjoint and their union is $V(K_n) \cup E(K_n)$.  (This gives yet another proof that $1+2+\cdots+n-1 = {{n}\choose{2}}$.)
 
Permute the elements of $K_n$ and then rank each team by the minimum position of its elements in the permutation. So for $K_3$, the permutation $13,1,12, 3, 2, 23$ induces the order: $T_3 < T_1 < T_2$.
{\it We claim that the probability $\cP$ that a random permutation will rank teams in the given order, 
$T_1 < T_2 < \cdots < T_n$,
is $\cP = 1/(2n-1)!!$}.

Suppose the claim is true.  Then the chance that one has a construction sequence with vertices in the given indexed order is the chance that the teams are in given order and further that in each team, the minimum position is occupied by the vertex.  
Thus, the above probability is multiplied by $1/n!$.  But the symmetry of $K_n$ means that each vertex sequence is as likely as any other, so one must multiply by $n!$, leaving the original probability of $1/(2n-1)!!$ as the probability that a randomly chosen permutation of the elements results in a construction sequence.

To verify the claim, $T_1$ has a probability of $2/n(n+1)$ of being first. Suppose that $T_1, T_2, \dots, T_{k-1}$ have been placed. Then there remain ${{n+1}\choose{2}} - {{k}\choose{2}} = \frac{(n+k)(n-k+1)}{2}$ unplaced elements, so the probability that some element in $T_k$ occupies the first available position is $2k/(n+k)(n-k+1)$.  But
\begin{equation}
\prod_{k=1}^n \frac{2k}{(n-k+1)(n+k)} = \frac{n! 2^n}{(2n)!} = \frac{1}{(2n-1)!!}
\label{eq:factorization}
\end{equation}
so the claim and the theorem are proved.
\end{proof}

Is the simple form of the probability that a random sequence is a c-sequence for $K_n$ just a consequence of equation (\ref{eq:factorization}) or is there a direct, non-algebraic argument?
For example, if $V(K_2)=\{a,b\}$ and if
 $y = (y_1,y_2,y_3) \in Perm(\{a,b,ab\})$, assuming uniform probability, $Pr(y_3 =ab)=\frac{1}{3}=\frac{1}{3!!}$.
Can this naive idea be extended?



\end{document}